\documentclass[11pt,a4paper]{article}

\title{\bf Asymptotic estimates on the time derivative of entropy on a Riemannian manifold}
\author{Adrian P. C. Lim$^a$,\quad Dejun Luo$^{a,b}$\footnote{Email: luodj@amss.ac.cn}
\vspace{3mm}\\
{\footnotesize $^a$UR Math\'{e}matiques, Universit\'{e} de
Luxembourg, 6, rue Richard Coudenhove-Kalergi, L-1359 Luxembourg}\\
{\footnotesize $^b$Key Lab of Random Complex Structures and Data
Science, Academy of Mathematics and}\\
{\footnotesize  Systems Science, Chinese Academy of Sciences,
Beijing 100190, China} }
\date{}

\usepackage{amssymb,amsmath,amsfonts,amsthm,color}

\def\H{\mathbb{H}}
\def\R{\mathbb{R}}

\def\d{\textup{d}}
\def\Ric{\textup{Ric}}
\def\Hess{\textup{Hess}}
\def\Ent{\textup{Ent}}
\def\Id{\textup{Id}}
\def\V{\textup{V}}

\def\ra{\rightarrow}
\def\<{\langle}
\def\>{\rangle}

\def\fin{\hfill$\square$}

\newtheorem{theorem}{Theorem}[section]
\newtheorem{lemma}[theorem]{Lemma}
\newtheorem{corollary}[theorem]{Corollary}
\newtheorem{proposition}[theorem]{Proposition}
\newtheorem{remark}[theorem]{Remark}

\begin{document}

\maketitle
\makeatletter 
\renewcommand\theequation{\thesection.\arabic{equation}}
\@addtoreset{equation}{section}
\makeatother 

\begin{abstract}
We consider the entropy of the solution to the heat equation on a
Riemannian manifold. When the manifold is compact, we provide two
estimates on the rate of change of the entropy in terms of the lower
bound on the Ricci curvature and the spectral gap respectively. Our explicit
computation for the three dimensional hyperbolic space shows that
the time derivative of the entropy is asymptotically bounded by two
positive constants.
\end{abstract}



\section{Introduction}

In Perelman's solution to the Poincar\'e conjecture, the
$W$-functional (i.e. the entropy of the heat kernel) played an
important role. Since then, there have been many attempts to
understand or generalize this functional to other situations, see
\cite{Chow, KotschwarNi, Ni04a, Ni04b}. This work is motivated by
\cite{Ni04a, Ni04b}, where the author presented the expression for
the time derivative of the entropy and studied its properties. Our
aim is to estimate the asymptotic behavior of the time derivative
of the entropy.

Let $M$ be a closed manifold, with dimension $n$. Equip $M$ with a
Riemannian metric $g$ and define the corresponding Laplacian $\Delta
\equiv\Delta_g$. We consider the following heat equation
  \begin{equation}\label{sect-1.1}
  \frac\partial{\partial t} u = \frac12\Delta u,\quad u_0=f,
  \end{equation}
where $f\in C^1(M,(0,\infty))$. In this work we follow the
probabilists' convention (e.g. \cite{ArnaudonThalmaier, Hsu02}) of
considering $\frac12\Delta$ instead of the Laplacian $\Delta$. Let
$\d x$ be the volume measure on $M$. Without loss of generality, we
assume that $M$ has volume 1. If $\int_M f\,\d x=1$, then $\int_M
u_t\,\d x\equiv 1$ for all $t>0$. Hence
  $$\int_M \frac{\partial}{\partial t}u_t\,\d x=\frac{\d}{\d t}\int_M
  u_t\,\d x=0.$$
Define the entropy
  \begin{equation*}
  \Ent(u_t)=-\int_M u_t\log u_t\,\d x.
  \end{equation*}
It is easy to see that $\Ent(u_t)\leq 0$ for all $t\geq0$, and the
equality holds if and only if $f\equiv 1$. By the integration by
parts formula, it is easy to show that
  \begin{align*}
  \frac{\d}{\d t}\Ent(u_t)&=\frac12\int_M\frac{|\nabla u_t|^2}{u_t}\,\d x.
  \end{align*}
From this, we see that if the initial value $f$ is not constant,
then $\Ent(u_t)$ is an increasing function of $t>0$. Denote by
$\Ric$ the Ricci curvature tensor on $M$. The main result of this
paper is

\begin{theorem}\label{sect-1-thm-1}
Assume that there is a $k\in\R$ such that $\Ric\geq k\Id$. Then
  \begin{equation}\label{sect-1-thm-1.1}
  \frac{\d}{\d t}\Ent(u_t)\leq
  \begin{cases}
  \frac{e^{-kt}}2 \Big[\frac1{q_0}-\frac{e^{-kt}-1}{nk}\Big]^{-1}, &
  k\neq0;\\
  \frac{nq_0}{2(n+q_0t)}, & k=0,
  \end{cases}
  \end{equation}
where $q_0=\int_{M}\frac{|\nabla f|^2}f\,\d x$.
\end{theorem}

This result will be proved in Section 2. We conclude from Theorem
\ref{sect-1-thm-1} that if the Ricci curvature is bounded from below
by a positive constant $k>0$, then the time derivative of the
entropy has an exponential decay as $t\ra\infty$; however, if $k<0$,
then the time derivative of the entropy is asymptotically dominated
by $-nk/2>0$ as $t\ra\infty$. The proof is based on the integration
by parts formula and the inequality \eqref{sect-2.5}. For a compact
Riemannian manifold with negative Ricci curvature, the estimate
given in Theorem \ref{sect-1-thm-1} is not satisfactory. In Section
3 we will give an estimate on the rate of change of entropy in terms
of the spectral gap, which implies that the time derivative of the
entropy always decrease exponentially to 0, even though Ricci is
bounded below by a negative constant (see Proposition \ref{p.l.1}).

One might wonder if a similar result holds when the manifold $M$ is
non-compact. The difficulty in this case is the justification of the
existence of the entropy and the integration by parts formula.
According to \cite[p.25]{KotschwarNi}, when the Ricci curvature of a
non-compact manifold is nonnegative, the entropy formula has been
established rigorously, see \cite{Chow}. Therefore by Li-Yau's
gradient estimate for positive solutions of the heat equation
\eqref{sect-1.1}, we can show that the time derivative of the
entropy does not exceed $n/(2t)$, i.e. the decay rate of the entropy
of heat kernel on $\R^n$; moreover, the time derivative achieves the
critical value $n/(2t)$ if and only if $M$ is isometric to $\R^n$,
see Theorem \ref{sect-2-thm-4}. However, when $M$ has negative Ricci
curvature, the situation is different. In the case of the three
dimensional hyperbolic space $\mathbb{H}^3$, we will prove in
Theorem \ref{sect-4-thm} that the rate of change of entropy is
asymptotically bounded by two positive constants.

\section{Rate of change of the entropy under curvature condition}

In this section we will provide the proof of Theorem
\ref{sect-1-thm-1}. In fact we can deal with more general cases. To
this end, let $Z$ be a $C^1$-vector field on $M$ and consider the
second order differential operator
  $$L=\frac12\Delta +Z.$$

\begin{proposition}\label{prop-1}
Let $f\in C^1(M,(0,\infty))$ and $u:[0,\infty)\times M \ra\R_+$ be a
solution to the heat equation
  \begin{equation}\label{prop-1.1}
  \frac\partial{\partial t} u = Lu,\quad u_0=f.
  \end{equation}
We have
  $$\big(L-\frac{\partial}{\partial t}\big)\bigg(\frac{|\nabla u|^2}{u}\bigg)
  =\frac{1}{u}\bigg|\Hess\,u -\frac{\nabla u \otimes \nabla u}{u} \bigg|^2 +
  \frac1u\big(\Ric(\nabla u, \nabla u)-2\< D_{\nabla u}Z,\nabla
  u\>\big).$$
\end{proposition}

\noindent{\bf Proof.} By the Weitzenb\"ock formula,
  \begin{equation}\label{prop-1.2}
  \frac12\Delta(|\nabla u|^2)=\<\nabla\Delta u,\nabla u\>
  +|\Hess\, u|_{HS}^2+\Ric(\nabla u,\nabla u),
  \end{equation}
where $\Hess\, u$ is the Hessian of $u$ and $|\cdot|_{HS}$ is the
Hilbert-Schmidt norm. Next taking $Y=\nabla u$ and $g=u$ in Lemma
\ref{lem-1} below gives us
  $$\<\nabla u,\nabla Z(u)\>=\< D_{\nabla u}Z,\nabla u\>+\<\nabla u,D_Z\nabla u\>,$$
where $D$ is the Levi-Civita connection. As a result,
  \begin{equation}\label{prop-1.3}
  Z(|\nabla u|^2)=Z\<\nabla u,\nabla u\>=2\<D_Z \nabla u,\nabla u\>
  =2\big\<\nabla Z(u)-D_{\nabla u}Z,\nabla u\big\>.
  \end{equation}
Finally $\frac\partial{\partial t}(|\nabla u|^2)=
2\big\<\frac\partial{\partial t}(\nabla u),\nabla u\big\>
=2\big\<\nabla(\frac\partial{\partial t} u),\nabla u\big\>$.
Combining this with \eqref{prop-1.2} and \eqref{prop-1.3}, and using
the equation \eqref{prop-1.1}, we get
  \begin{equation}\label{prop-1.4}
  \bigg(L-\frac\partial{\partial t}\bigg)(|\nabla u|^2)
  =|\Hess\,u|_{HS}^2+\Ric(\nabla u,\nabla u)-2\< D_{\nabla u}Z,\nabla
  u\>.
  \end{equation}

Now notice that
  \begin{equation}\label{prop-1.5}
  \Delta\bigg(\frac{|\nabla u|^2}u\bigg)=u^{-1}\Delta(|\nabla u|^2)
  +|\nabla u|^2\Delta(u^{-1})+2\big\<\nabla u^{-1},\nabla(|\nabla
  u|^2)\big\>.
  \end{equation}
We have
  $$\Delta(u^{-1})=-\frac{\Delta u}{u^2}+\frac{2}{u^3}|\nabla u|^2$$
and
  $$\big\<\nabla u^{-1},\nabla(|\nabla u|^2)\big\>
  =-\frac1{u^2}\<\nabla u,\nabla(|\nabla u|^2)\>
  =-\frac2{u^2}\< D_{\nabla u}\nabla u,\nabla u\>
  =-\frac2{u^2}(\Hess\,u)(\nabla u,\nabla u).$$
Substituting the above equalities into \eqref{prop-1.5} gives rise
to
  \begin{equation}\label{prop-1.6}
  \Delta\bigg(\frac{|\nabla u|^2}u\bigg)=\frac1u\Delta(|\nabla u|^2)
  -\frac{|\nabla u|^2}{u^2}\Delta u+\frac2{u^3}|\nabla u|^4
  -\frac4{u^2}(\Hess\,u)(\nabla u,\nabla u).
  \end{equation}
We also have
  \begin{align}\label{prop-1.7}
  Z\bigg(\frac{|\nabla u|^2}u\bigg)&=\big\< Z,|\nabla u|^2\nabla(u^{-1})+u^{-1}\nabla(|\nabla
  u|^2)\big\>\cr
  &=-\frac{|\nabla u|^2}{u^2}Z(u)+\frac1u Z(|\nabla u|^2)
  \end{align}
and
  \begin{equation}\label{prop-1.8}
  \frac\partial{\partial t}\bigg(\frac{|\nabla u|^2}u\bigg)
  =-\frac{|\nabla u|^2}{u^2}\frac\partial{\partial t}u
  +\frac1u \frac\partial{\partial t}(|\nabla u|^2).
  \end{equation}
Consequently, by \eqref{prop-1.6}--\eqref{prop-1.8},
  \begin{align*}
  \bigg(L-\frac{\partial}{\partial t}\bigg)\bigg(\frac{|\nabla u|^2}{u}\bigg)
  &=\frac1u\bigg(L-\frac{\partial}{\partial t}\bigg)(|\nabla u|^2)
  -\frac{|\nabla u|^2}{u^2}\bigg(L-\frac{\partial}{\partial t}\bigg)u\cr
  &\hskip12pt+\frac{|\nabla u|^4}{u^3}-\frac2{u^2}(\Hess\,u)(\nabla u,\nabla
  u),
  \end{align*}
which, by \eqref{prop-1.4}, is equal to
  \begin{align*}
  &\hskip12pt \frac1u\big(|\Hess\,u|_{HS}^2+\Ric(\nabla u,\nabla u)-2\< D_{\nabla u}Z,\nabla
  u\>\big)+\frac{|\nabla u|^4}{u^3}-\frac2{u^2}(\Hess\,u)(\nabla u,\nabla
  u)\cr
  &= \frac1u\bigg(\Big|\Hess\,u-\frac{\nabla u\otimes\nabla u}u\Big|_{HS}^2
  +\Ric(\nabla u,\nabla u)-2\< D_{\nabla u}Z,\nabla u\>\bigg).
  \end{align*}
The proof is complete. \fin

\begin{lemma}\label{lem-1}
Let $Y,Z$ be two $C^1$-vector fields on $M$ and $g\in C^2(M)$. We
have
  $$\< Y,\nabla Z(g)\>=\< D_Y Z,\nabla g\>+\< Y,D_Z\nabla g\>.$$
\end{lemma}

\noindent{\bf Proof.} By the consistency of the Levi-Civita
connection $D$ with the metric,
  $$\< Y,\nabla Z(g)\>=Y[Z(g)]=Y\< Z,\nabla g\>=\< D_Y Z,\nabla g\>
  +\< Z,D_Y \nabla g\>.$$
Notice that
  $$\< Z,D_Y \nabla g\>=(\Hess\,g)(Y,Z)=(\Hess\,g)(Z,Y)=\< Y,D_Z \nabla g\>.$$
The result follows. \fin

\medskip

Now we take $V\in C^2(M)$ and let $Z=\nabla V$. That is, we consider
the heat equation
  \begin{equation}\label{sect-2.1}
  \frac\partial{\partial t}u=Lu=\frac12\Delta u+\<\nabla V,\nabla u\>,
  \quad u_0=f>0.
  \end{equation}
Define the measure $\d\mu=e^{2V}\,\d x$ on $M$, where $\d x$ is the
volume element of $M$. Normalize the measure $\mu$ if necessary,
we will assume that $\mu(M)=1$. For any $f,g\in C^1(M)$, the following
integration by parts formula holds:
  \begin{equation}\label{sect-2.2}
  -\int_M Lf\cdot g\,\d\mu=\frac12\int_M\<\nabla f,\nabla
  g\>\,\d\mu.
  \end{equation}
Assume that $\mu(f)=\int_M f\,\d\mu=1$. Then we have $\int_M
u_t\,\d\mu\equiv 1$ for all $t>0$. Hence
  $$\int_M \frac{\partial}{\partial t}u_t\,\d\mu=\frac{\d}{\d t}\int_M
  u_t\,\d\mu=0.$$
Define the entropy
  \begin{equation}\label{sect-2.3}
  \Ent(u_t)=-\int_M u_t\log u_t\,\d\mu.
  \end{equation}
It is easy to see that $\Ent(u_t)\leq 0$ for all $t\geq0$, and the
equality holds if and only if $f\equiv 1$. Using the integration by
parts formula \eqref{sect-2.2},
  \begin{align}\label{sect-2.4}
  \frac{\d}{\d t}\Ent(u_t)&=-\int_M\bigg((\log u_t)\frac{\partial}{\partial t}
  u_t+\frac{\partial}{\partial t}u_t\bigg)\d\mu=-\int_M (\log
  u_t)Lu_t\,\d\mu\cr
  &=\frac12\int_M\<\nabla\log u_t,\nabla u_t\>\,\d\mu
  =\frac12\int_M\frac{|\nabla u_t|^2}{u_t}\d\mu.
  \end{align}
From this, we see that if the initial value $f$ is not a constant,
then $\Ent(u_t)$ is an increasing function of $t>0$.

We want to estimate the rate of change of the entropy, i.e.
$\frac{\d}{\d t}\Ent(u_t)$, as $t\ra\infty$. To this end, we will
make use of the equality proved in Proposition \ref{prop-1}. It is
easy to see that
  \begin{equation}\label{sect-2.4.5}
  \bigg|\Hess\,u-\frac{\nabla u\otimes\nabla u}u\bigg|_{HS}^2
  =u^2|\Hess(\log u)|_{HS}^2\geq \frac{u^2}n|\Delta(\log u)|^2,
  \end{equation}
where $n$ is the dimension of $M$, and the equality holds if and
only if $\Hess(\log u)=g\,\Id$ with some function $g:\R_+\times
M\ra\R$. Therefore by Proposition \ref{prop-1}, we arrive at
  \begin{equation}\label{sect-2.5}
  \bigg(L-\frac{\partial}{\partial t}\bigg)\bigg(\frac{|\nabla u|^2}{u}\bigg)
  \geq \frac un|\Delta(\log u)|^2+\frac1u\big(\Ric-2\,\Hess\,V\big)(\nabla u,\nabla
  u).
  \end{equation}

\begin{theorem}\label{sect-2-thm-1}
Assume that there is a $k>0$ such that
  $$\Ric-2\,\Hess\,V\geq k\,\Id.$$
Then for all $t>0$,
  $$\frac{\d}{\d t}\Ent(u_t)\leq\frac{ e^{-kt}}2\int_{M}\frac{|\nabla f|^2}f\,\d\mu.$$
\end{theorem}

\noindent{\bf Proof.} Integrating both sides of \eqref{sect-2.5}, we
obtain
  $$\int_M\bigg(L-\frac{\partial}{\partial t}\bigg)\bigg(\frac{|\nabla u|^2}{u}\bigg)\d\mu
  \geq\frac 1n\int_M u|\Delta(\log u)|^2\,\d\mu+k\int_M\frac{|\nabla u|^2}u\,\d\mu.$$
By \eqref{sect-2.2}, we have
  $$\int_M L\bigg(\frac{|\nabla u|^2}{u}\bigg)\d\mu=0,$$
thus the above inequality reduces to
  \begin{equation}\label{sect-2-thm-1.1}
  -\int_M\frac{\partial}{\partial t}\bigg(\frac{|\nabla u|^2}{u}\bigg)\d\mu
  \geq\frac 1n\int_M u|\Delta(\log u)|^2\,\d\mu+k\int_M\frac{|\nabla u|^2}u\,\d\mu.
  \end{equation}
Therefore
  $$\frac{\d}{\d t}\int_M\frac{|\nabla u|^2}{u}\,\d\mu\leq
  -k\int_M\frac{|\nabla u|^2}{u}\,\d\mu,$$
from which the desired estimate follows. \fin

\medskip

Now we can prove Theorem \ref{sect-1-thm-1}, the case when $V = 0$.

\medskip

\noindent{\bf Proof of Theorem \ref{sect-1-thm-1}.} Remark that now
$\d\mu=\d x$, i.e. the volume measure. We still have
\eqref{sect-2-thm-1.1}, or
  \begin{equation*}
  -\frac{\d}{\d t}\int_M\frac{|\nabla u|^2}{u}\,\d x
  \geq\frac 1n\int_M u|\Delta(\log u)|^2\,\d x+k\int_M\frac{|\nabla u|^2}u\,\d x.
  \end{equation*}
For simplicity of notation, we denote by $q_t=\int_M \frac{|\nabla
u_t|^2}{u_t}\,\d x$. Then the above inequality can be written as
  \begin{equation}\label{sect-2-thm-2.3}
  -\frac{\d q_t}{\d t}\geq\frac 1n\int_M u|\Delta(\log u)|^2\,\d x+kq_t.
  \end{equation}
By the Cauchy inequality,
  \begin{equation}\label{sect-2-thm-2.3.5}
  \bigg(\int_M u\,\Delta \log u\,\d x\bigg)^2
  \leq\int_M u|\Delta \log u|^2\,\d x \cdot \int_M u\,\d x
  =\int_M u|\Delta \log u|^2\,\d x.
  \end{equation}
Using the integration by parts formula, we have
  $$\int_M u\,\Delta \log u\,\d x=-\int_M\<\nabla u,\nabla\log u\>\,\d x
  =-q_t.$$
Therefore \eqref{sect-2-thm-2.3} becomes
  \begin{equation}\label{sect-2-thm-2.4}
  \frac{\d q_t}{\d t}\leq -\frac{q_t^2}n-kq_t.
  \end{equation}

Now we solve the above differential inequality. If $k=0$, then
  $$\frac{\d q_t}{\d t}\leq -\frac{q_t^2}n,$$
from which we easily obtain that $q_t\leq \frac{nq_0}{n+q_0t}$. In
the case $k\neq0$, we divide both sides of \eqref{sect-2-thm-2.4} by
$q_t^2$ and get
  $$\frac1{q_t^2}\frac{\d q_t}{\d t} \leq -\frac1n-\frac k{q_t}.$$
This is equivalent to
  $$\frac{\d q_t^{-1}}{\d t}\geq \frac1n+kq_t^{-1}.$$
Therefore
  $$\frac{\d}{\d t}\big(e^{-kt}q_t^{-1}\big)\geq \frac{e^{-kt}}n.$$
Integrating this inequality from 0 to $t$ leads to
  $$e^{-kt}q_t^{-1}-q_0^{-1}\geq \frac1{nk}(1-e^{-kt}).$$
The proof is now completed. \fin

\medskip

We have the following simple observations.

\begin{corollary} Assume the conditions of Theorem
\ref{sect-1-thm-1}.
\begin{itemize}
\item[\rm(i)] If $k>0$, then the time derivative of the entropy has
an exponential decay as $t\ra\infty$;
\item[\rm(ii)] if $k<0$, then the time derivative of the entropy
is asymptotically dominated by $-nk/2$ as $t\ra\infty$.
\end{itemize}
\end{corollary}

\noindent{\bf Proof.} The assertions follow directly from Theorem
\ref{sect-1-thm-1}. \fin

\medskip

We can also obtain an estimate on the rate of change of entropy from
Hamilton's gradient estimate for heat equations (see
\cite{Hamilton94} or \cite[Corollary 3.3]{ArnaudonThalmaier}).

\begin{proposition}\label{sect-2-prop-1}
Let $u$ be a solution to the standard heat equation
\eqref{sect-1.1}. Assume $\int_M f\,\d x=1$ and there is $k\in\R$
such that $\Ric\geq k\,\Id$. Then for all $t>0$,
  $$\frac{\d}{\d t}\Ent(u_t)\leq \bigg(\frac1t -k\bigg)\log(\sup f).$$
Moreover, if $\,\Ric\geq0$, then we can take $k=0$ in the above
estimate.
\end{proposition}

\noindent{\bf Proof.} By \cite[Corollary 3.3]{ArnaudonThalmaier}, we
have
  $$\frac{|\nabla u_t|^2}{u_t^2}\leq 2\bigg(\frac1t -k\bigg)\log\frac A{u_t},$$
where $A=\sup_{M\times [0,t]}u=\sup_M f$. Therefore
  $$\frac{|\nabla u_t|^2}{u_t}\leq 2\bigg(\frac1t -k\bigg)u_t\log\frac{\sup f}{u_t}.$$
Integrating both sides on $M$, we get
  \begin{align}
  \int_M\frac{|\nabla u_t|^2}{u_t}\,\d x&\leq 2\bigg(\frac1t -k\bigg)
  \int_M u_t\log\frac{\sup f}{u_t}\,\d x\cr
  &=2\bigg(\frac1t -k\bigg)\big(\log({\sup f})+\Ent(u_t)\big)
  \leq 2\bigg(\frac1t -k\bigg)\log({\sup f}),
  \end{align}
since the entropy $\Ent(u_t)\leq 0$. The proof is complete. \fin

\medskip

It is interesting to compare the two estimates given in Theorem
\ref{sect-1-thm-1} and Proposition \ref{sect-2-prop-1}. We
distinguish three cases:
\begin{itemize}
\item[\rm(a)] $k<0$. By Theorem \ref{sect-1-thm-1}, $\frac{\d}{\d
t}\Ent(u_t)$ is asymptotically dominated by $-nk/2$, depending only
on the dimension and curvature; while by Proposition
\ref{sect-2-prop-1}, the asymptotic constant is $-k\log(\sup f)$
which depends on the initial condition.

\item[\rm(b)] $k=0$. By Theorem \ref{sect-1-thm-1}, $\frac{\d}{\d
t}\Ent(u_t)\leq \frac{nq_0}{2(n+q_0t)}\leq \frac n{2t}$. The
estimate given by Proposition \ref{sect-2-prop-1} is $\frac{\d}{\d
t}\Ent(u_t)\leq \frac1t\log(\sup f)$.

\item[\rm(c)] $k>0$. Theorem \ref{sect-1-thm-1} gives us an
exponential decay: $\frac{\d}{\d t}\Ent(u_t)\leq \frac12 q_0
e^{-kt}$; while Proposition \ref{sect-2-prop-1} only leads to a
polynomial decay: $\frac{\d}{\d t}\Ent(u_t)\leq \frac1t\log(\sup
f)$.
\end{itemize}

To complete this section, we briefly discuss the entropy of the heat
kernel on a non-compact Riemannian manifold with nonnegative Ricci
curvature. As mentioned in the Introduction (see also
\cite[p.90]{Ni04a}), the integration by parts formula can be
justified rigorously, thanks to Li-Yau's gradient estimate
\cite{LiYau}. Therefore if $u_t$ is the heat kernel of $M$, in the
same way we can show that $\frac{\d}{\d t}\Ent(u_t)\leq \frac
n{2t}$. The next theorem is an analogue of \cite[Theorem
1.4]{Ni04a}, which says that if the Ricci curvature is nonnegative,
then the decay rate of $\frac n{2t}$ is achieved if and only if
$M=\R^n$.

\begin{theorem}\label{sect-2-thm-4}
Let $M$ be a non-compact Riemannian manifold with nonnegative Ricci
curvature, and $u_t$ its heat kernel. Then
  \begin{align}\label{sect-2-thm-3}
  \frac{\d}{\d t}\Ent(u_t)=\frac{n}{2t}
  \end{align}
if and only if $M$ is isometric to $\R^n$.
\end{theorem}

\noindent{\bf Proof.} We only have to prove the necessity part.
Suppose \eqref{sect-2-thm-3} holds. Since $M$ has nonnegative Ricci
curvature, it is stochastically complete, that is, $\int_M u_t\,\d
x=1$ for all $t>0$ (see e.g. \cite[Theorem 4.2.4]{Hsu02}). By the
proof of Theorem \ref{sect-1-thm-1} (remark that now $k=0$), we must
have equalities in \eqref{sect-2.4.5} and \eqref{sect-2-thm-2.3.5}.
The equality in \eqref{sect-2.4.5} will imply that $\Delta(\log
u_t)= ng_t$ for some function $g:\R_+\times M \ra\R$. Next the
equality in \eqref{sect-2-thm-2.3.5} means that $g_t$ is a function
independent of $x\in M$. By the integration by parts formula,
  $$ng_t=\int_M u_t\,\Delta(\log u_t)\,\d x=-\int_M\frac{|\nabla u_t|^2}{u_t}\,\d x
  =-2\frac{\d}{\d t}\Ent(u_t)=-\frac{n}{t},$$
therefore $g_t=-\frac1t$ and $\Delta(\log u_t)=-\frac nt$.

The rest of the proof is similar to that of \cite[Theorem
1.4]{Ni04a}. By Varadhan's large deviation formula (see
\cite[Theorem 5.2.1]{Hsu02}),
  \begin{align*}
  -2\lim_{t \rightarrow 0} t\log u_t(x,y) = d^2(x,y),
  \end{align*}
where $d:M\times M\ra\R_+$ is the Riemannian distance function.
Thus,
  \begin{align}\label{sect-2-thm-3.1}
  \Delta d^2(x,y) = -2\lim_{t \ra0} t\,\Delta\log u_t(x,y) = 2n.
  \end{align}
From \eqref{sect-2-thm-3.1} we deduce that
  $$\frac{A_x(r)}{V_x(r)}=n,$$
where $A_x(r)$ and $V_x(r)$ denote respectively the area of
$\partial B_x(r)$ and the volume of $B_x(r)=\{y\in M:d(x,y)\leq
r\}$. This implies that $V_x(r)$ is the same as the volume function
of Euclidean balls. The equality case of the volume comparison
theorem gives us $M=\R^n$. \fin

\section{Rate of change of the entropy in terms of the spectral gap}

In this section we will obtain an estimate on the time derivative of
the entropy in terms of the spectral gap. For simplicity, set $V = 0$.
Then we have the Sturm-Liouville decomposition of the fundamental
solution,
  \begin{equation}\label{sect-3.1}
  \Phi(t,x,y) = \sum_{j=0}^\infty e^{-\lambda_j t/2}\phi_j(x)\phi_j(y)
  \end{equation}
whereby $\phi_j$ are the eigenfunctions of $\Delta$ with eigenvalue
$\lambda_j$, $0 = \lambda_0 < \lambda_1 < \ldots $, $\lambda_j
\nearrow \infty$. In particular, each eigenvalue has finite
multiplicity and each $\phi_j$ is smooth. Furthermore,
$\{\phi_j\}_{j=0}^\infty$ form an orthonormal basis for the $L^2(M)$
and $\phi_0 = 1$ if $\int_M \Phi_t\, \d y = 1$, see
\cite[p.139]{Chavel84}. The series in \eqref{sect-3.1}
converges absolutely and uniformly.

Given any solution $u$ to the heat equation with initial data $f$, we can
write
  \begin{equation}\label{sect-3.2}
  u_t(x) = \int_M f(y)\Phi(t,x,y)\, \d y = \sum_{j=0}^\infty e^{-\lambda_j t/2}c_j\phi_j(x),\ c_j
  = \int_M f(y)\phi_j(y)\, \d y.
  \end{equation}
And
  \begin{equation}\label{sect-3.3}
  \frac\partial{\partial t} u_t(x) = \frac1{2}\Delta u_t(x) =
  -\sum_{j=0}^\infty e^{-\lambda_j t/2}\lambda_jc_j\phi_j(x)/2.
  \end{equation}
The series in \eqref{sect-3.2} and \eqref{sect-3.3}
converge absolutely and uniformly for $t > 0$.

Write $\langle f, g \rangle = \int_M fg\, \d x$ and $\| f \|_2^2 =
\langle f, f\rangle$.

\begin{proposition}\label{p.l.1}
Assume $M$ is a closed Riemannian manifold. Let $f \in
C^2(M,(0,\infty))$ and $u_t$ be the solution to the heat equation
with $u_0 = f$. Then
  $$ 0 \leq \frac{\d}{\d t}\Ent(u_t) \leq \frac12 e^{-\lambda_1 t/2}\| \Delta f \|_2
  \sqrt{\textup{Vol}(M)} \,(|\log\inf  f| + |\log\sup  f|). $$
Here, $\textup{Vol}(M)$ is the volume of the manifold.
\end{proposition}

\noindent{\bf Proof.} By our assumption, $\Delta f$ is continuous.
Thus,
  \begin{equation*}
  \Delta f = \sum_{j=0}^\infty \langle \Delta f,
  \phi_j \rangle \phi_j = -\sum_{j=0}^\infty \lambda_j\langle f,
  \phi_j \rangle\phi_j = -\sum_{j=1}\lambda_j c_j \phi_j,
  \end{equation*}
and $\|\Delta f\|_2^2 = \langle \Delta f,\Delta f \rangle  =
\sum_{j=1} \lambda_j^2 c_j^2 < \infty$. By \eqref{sect-3.3},
  $$\|\Delta u_t\|_2^2=\sum_{j=1}^\infty e^{-\lambda_j t}\lambda_j^2 c_j^2
  \leq e^{-\lambda_1 t}\sum_{j=1}^\infty\lambda_j^2 c_j^2
  =e^{-\lambda_1 t}\|\Delta f\|_2^2.$$
Next by \eqref{sect-3.2}, it is clear that
  \begin{align*}
  \int_M |\log u|^2\, \d x \leq \textup{Vol}(M) \left(|\log\inf  f| + |\log\sup  f|
  \right)^2.
  \end{align*}
Now the estimate follows from Cauchy's inequality:
  \begin{align*}
  \frac{\d}{\d t}\Ent(u_t)
  &= -\frac12\int_M \log u_t \cdot \Delta u_t\, \d\nu \cr
  &\leq \frac12\|\log u_t\|_2\|\Delta u_t\|_2\cr
  &\leq \frac12e^{-\lambda_1 t/2}\| \Delta f \|_2
  \sqrt{\textup{Vol}(M)} \,(|\log\inf  f| + |\log\sup  f|).
  \end{align*}

\section{Example: Entropy of the heat kernel on the three dimensional
hyperbolic space $\H^3$}

In this section we consider the three dimensional hyperbolic space
$\H^3$ of constant sectional curvature $k<0$. Let $d(\cdot, \cdot)$
be the Riemannian distance function on $\H^3$. The heat kernel on
$\H^3$ has the following explicit formula (see
\cite[p.150]{Chavel84}):
  $$h(t,x,y)=e^{-d(x,y)^2/2t}(2\pi t)^{-3/2}\frac{\sqrt{-k}\, d(x,y)}{\sinh\sqrt{-k}\, d(x,y)}
  e^{k t/2},\quad x,y\in\H^3.$$
(There is a mistake in the formula for $h(t,x,y)$ given in
\cite[Example 5.1.3]{Hsu02}: there the last factor is
$e^{-t}\,(k=-1)$, rather than $e^{-t/2}$.) Let $\V_{\H^3}$ be the
volume measure on $\H^3$ and define the entropy
  $$\Ent(h(t,x,\cdot))=-\int_{\H^3}h(t,x,y)\log h(t,x,y)\,\d\V_{\H^3}(y),
  \quad (t,x)\in (0,\infty)\times \H^3.$$
The main result in this section is

\begin{theorem}\label{sect-4-thm}
  \begin{align*}
  -k(2-\log\sqrt 2\,)\leq \varliminf_{t\ra\infty}\frac{\d}{\d t}\Ent(h(t,x,\cdot))
  \leq\varlimsup_{t\ra\infty}\frac{\d}{\d t}\Ent(h(t,x,\cdot))\leq
  -k(2+\log\sqrt 2\,).
  \end{align*}
\end{theorem}

To prove this theorem, we need some preparations. We denote by
$\kappa=\sqrt{-k}\,$ for the simplification of notations. Using
polar coordinates the metric on $\H^3$ is expressed as
  $$\d s^2=\d r^2+\kappa^{-2}(\sinh \kappa r)^2\d\theta^2,$$
where $\d\theta$ is the standard volume measure on the sphere $S^2$.
Hence for any integrable function $f:\H^3\ra\R$, the following
equality holds:
  \begin{equation}\label{sect-4.1}
  \int_{\H^3}f(x)\,\d\V_{\H^3}(x)=\int_0^\infty\bigg(\int_{S^2}f(r\theta)\,\d\theta\bigg)
  \frac{(\sinh\kappa r)^2}{\kappa^2}\,\d r.
  \end{equation}
We have
  \begin{align}\label{sect-4.2}
  \Ent(h(t,x,\cdot))=\frac32\log(2\pi t)+\frac{\kappa^2t}2+I_1(t)+I_2(t),
  \end{align}
where
  \begin{align}
  I_1(t)&=\frac1{2t}\int_{\H^3}h(t,x,y)d(x,y)^2\,\d\V_{\H^3}(y),\label{sect-4.3}\\
  I_2(t)&=\int_{\H^3}h(t,x,y)\log\frac{\sinh\kappa
  d(x,y)}{\kappa d(x,y)}\,\d\V_{\H^3}(y).\label{sect-4.4}
  \end{align}

We first compute $I_1(t)$. Using the expression of $h(t,x,y)$ and by
the formula \eqref{sect-4.1}, we have
  \begin{align*}
  I_1(t)&=\frac1{2t(2\pi t)^{3/2}e^{\kappa^2t/2}}\int_0^\infty
  \bigg(\int_{S^2}e^{-r^2/2t}\frac{\kappa r}{\sinh\kappa r}\,r^2\,\d\theta\bigg)
  \frac{(\sinh\kappa r)^2}{\kappa^2}\,\d r\cr
  &=\frac1{\sqrt{2\pi}\,\kappa t^{5/2}e^{\kappa^2t/2}}\,I_{11}(t),
  \end{align*}
where $I_{11}(t)=\int_0^\infty e^{-r^2/2t}r^3\sinh\kappa r\,\d r$.
Here we collect some results for later use.

\begin{lemma}\label{sect-4-lem-1}
Let $\alpha(t)=\int_0^{\kappa t^{1/2}}e^{-r^2/2}\,\d r$. Then
  \begin{eqnarray*}
  \int_0^\infty e^{-r^2/2t}\sinh\kappa r\,\d r &=& t^{1/2}e^{\kappa^2 t/2}
  \alpha(t),\cr
  \int_0^\infty e^{-r^2/2t}\cosh\kappa r\,\d r &=& \sqrt{\frac\pi2}\,t^{1/2}e^{\kappa^2 t/2}  ,\cr
  \int_0^\infty e^{-r^2/2t}r\sinh\kappa r\,\d r &=& \sqrt{\frac\pi2}\,\kappa t^{3/2}e^{\kappa^2t/2},\cr
  \int_0^\infty e^{-r^2/2t}r\cosh\kappa r\,\d r &=& t+\kappa t^{3/2}e^{\kappa^2 t/2}\alpha(t) ,\cr
  \int_0^\infty e^{-r^2/2t}r^2\sinh\kappa r\,\d r &=& \kappa t^2+t^{3/2}(\kappa^2t+1)e^{\kappa^2 t/2}\alpha(t),\cr
  \int_0^\infty e^{-r^2/2t}r^2\cosh\kappa r\,\d
  r &=&\sqrt{\frac{\pi}2}\,t^{3/2}(\kappa^2t+1)e^{\kappa^2 t/2},\cr
  \int_0^\infty e^{-r^2/2t}r^3\sinh\kappa r\,\d r &=&\sqrt{\frac{\pi}2}
  \,\kappa t^{5/2}(\kappa^2 t+3)e^{\kappa^2 t/2},\cr
  \int_0^\infty e^{-r^2/2t}r^3\cosh\kappa r\,\d r &=& t^2(\kappa^2 t+2)
  +\kappa t^{5/2}(\kappa^2 t+3)e^{\kappa^2 t/2}\alpha(t),\cr
  \int_0^\infty e^{-r^2/2t}r^4\sinh\kappa r\,\d r &=&\kappa t^3(\kappa^2 t+5)
  +t^{5/2}(\kappa^4 t^2+6\kappa^2 t+3)e^{\kappa^2 t/2}\alpha(t).
  \end{eqnarray*}
\end{lemma}

\noindent{\bf Proof.} We only prove the first two equalities. The
others can be proved using the integration by parts formula. We
have
  $$\int_0^\infty e^{-r^2/2t}e^{\kappa r}\,\d r=e^{\kappa^2 t/2}\int_{-\kappa t}^\infty e^{-r^2/2t}\,\d r
  =t^{1/2}e^{\kappa^2 t/2}\int_{-\kappa t^{1/2}}^\infty e^{-r^2/2}\,\d r$$
and
  $$\int_0^\infty e^{-r^2/2t}e^{-\kappa r}\,\d r=e^{\kappa^2 t/2}\int_{\kappa t}^\infty e^{-r^2/2t}\,\d r
  =t^{1/2}e^{\kappa^2 t/2}\int_{\kappa t^{1/2}}^\infty e^{-r^2/2}\,\d r.$$
Hence
  $$\int_0^\infty e^{-r^2/2t}\sinh r\,\d r=\frac12t^{1/2}e^{\kappa^2 t/2}
  \int_{-\kappa t^{1/2}}^{\kappa t^{1/2}} e^{-r^2/2}\,\d r=t^{1/2}e^{\kappa^2 t/2}
  \alpha(t).$$

Next since the function $r\mapsto\cosh\kappa r$ is even,
  \begin{align*}
  \int_0^\infty e^{-r^2/2t}\cosh\kappa r\,\d r
  &=\frac12\int_{-\infty}^\infty e^{-r^2/2t}\cosh r\,\d r\cr
  &=\frac12\int_{-\infty}^\infty e^{-r^2/2t}e^{\kappa r}\,\d r
  =\sqrt{\frac\pi2}\,t^{1/2}e^{\kappa^2 t/2}.
  \end{align*}

\fin

\medskip

By Lemma \ref{sect-4-lem-1}, we obtain
  \begin{equation}\label{sect-4.5}
  I_1(t)=\frac1{\sqrt{2\pi}\,\kappa t^{5/2}e^{\kappa^2t/2}}\,\sqrt{\frac{\pi}2}
  \,\kappa t^{5/2}(\kappa^2 t+3)e^{\kappa^2 t/2} =\frac12 (\kappa^2 t+3).
  \end{equation}

Now we consider $I_2(t)$. Again by \eqref{sect-4.1},
  \begin{equation}\label{sect-4.6}
  I_2(t)=\sqrt{\frac2{\pi}}\,\kappa^{-1} t^{-3/2}e^{-\kappa^2 t/2}\int_0^\infty
  e^{-2r^2/2t}r\sinh\kappa r \log\frac{\sinh\kappa r}{\kappa r}\,\d r.
  \end{equation}
Due to the presence of the term $\log\frac{\sinh\kappa r}{\kappa
r}$, we are unable to compute $I_2(t)$ explicitly. In the sequel we
intend to find some estimates on it. Define
  \begin{equation}\label{sect-4.7}
  \xi(t)=\sqrt{\frac2{\pi}}\,\kappa^{-1} t^{-3/2}e^{-\kappa^2 t/2}\quad \mbox{and} \quad
  \eta(t)=\int_0^\infty e^{-r^2/2t}r\sinh\kappa r \log\frac{\sinh\kappa r}{\kappa r}\,\d r.
  \end{equation}
We have
  \begin{equation}\label{sect-4.8}
  \xi^\prime(t)=-\frac1{\sqrt{2\pi}\,\kappa}\cdot\frac{\kappa^2 t+3}{t^{5/2}e^{\kappa^2 t/2}}<0.
  \end{equation}
Therefore, to estimate $I_2^\prime(t)$, it is enough to estimate
$\eta(t)$ and $\eta^\prime(t)$.

The following lemma gives the key ingredient.

\begin{lemma}\label{sect-4-lem-2}
For any $r>0$,
  $$\frac1{1+2r}<\frac{1-e^{-2r}}{2r}<\frac1{1+r}.$$
\end{lemma}

\noindent{\bf Proof.} (1) Let $\phi(s)=1-e^{-s}-se^{-s}$. Then
$\phi(0)=0$ and $\phi^\prime(s)=se^{-s}>0$. Hence $\phi(s)>0$ for
all $s>0$. This implies
  $$1+s-e^{-s}-se^{-s}>s,$$
which leads to the first inequality by taking $s=2r$.

(2) Let $\psi(r)=1-r-e^{-2r}-re^{-2r}$. Then $\psi(0)=0$ and
$\psi^\prime(r)=-1+e^{-2r}+2re^{-2r}$. We have $\psi^\prime(0) =0$
and $\psi''(r)=-4re^{-2r}<0$ for any $r>0$. Therefore
$\psi^\prime(r)<0$. As a result, $\psi(r)<0$ for all $r>0$. This
implies that
  $$1+r-e^{-2r}-re^{-2r}<2r,$$
which is equivalent to the second inequality. \fin

\begin{remark}\label{sect-4-rem-1}
Fix any $\beta\in(1,2)$. In the same way we can show that for
$r\in\big(0,\frac{\beta-1}\beta\big)$, it holds
$\frac{1-e^{-2r}}{2r}> \frac1{1+\beta r}$. Moreover, it is clear
that $\frac{1-e^{-2r}}{2r}<\frac1{1+\beta r}$ when $r$ is
sufficiently large. Therefore the two inequalities in Lemma
\ref{sect-4-lem-2} are sharp in this sense.
\end{remark}

Notice that
  $$\frac{\sinh\kappa r}{\kappa r}=\frac{e^{\kappa r}(1-e^{-2\kappa r})}{2\kappa r},$$
thus by Lemma \ref{sect-4-lem-2},
  \begin{equation}\label{sect-4.9}
  \kappa r+\log\frac1{1+2\kappa r}<\log\frac{\sinh\kappa r}{\kappa r}
  <\kappa r+\log\frac1{1+\kappa r},\quad r>0.
  \end{equation}

Now we can obtain the upper and lower bound on $\eta(t)$.

\begin{lemma}\label{sect-4-lem-3}
  \begin{align*}
  \eta(t)&> \kappa^2t^2+\kappa t^{3/2}(\kappa^2t+1)e^{\kappa^2t/2}\alpha(t)
  -\sqrt{\frac{\pi}2}\,\kappa t^{3/2}e^{\kappa^2 t/2}\log(2\kappa^2t+4),\cr
  \eta(t)&< \kappa^2t^2+\kappa t^{3/2}(\kappa^2t+1)e^{\kappa^2t/2}\alpha(t)
  -\sqrt{\frac{\pi}2}\,\kappa t^{3/2}e^{\kappa^2t/2}
  \log\bigg(1+\sqrt{\frac{\pi}2}\,\kappa^2t\alpha(t)^{-1}\bigg).
  \end{align*}
\end{lemma}

\noindent{\bf Proof.} By the definition \eqref{sect-4.7} and the
inequality \eqref{sect-4.9}, we have
  \begin{equation}\label{sect-4-lem-3.1}
  \eta(t)>\eta_1(t)-\eta_2(t),
  \end{equation}
where by Lemma \ref{sect-4-lem-1},
  \begin{equation}\label{sect-4-lem-3.2}
  \eta_1(t)=\kappa\int_0^\infty e^{-r^2/2t}r^2\sinh\kappa r\,\d r
  =\kappa^2t^2+\kappa t^{3/2}(\kappa^2t+1)e^{\kappa^2t/2}\alpha(t)
  \end{equation}
and $\eta_2(t)=\int_0^\infty e^{-r^2/2t}r\log(1+2\kappa r)\sinh
\kappa r\,\d r$. Define the measure
  $$\d\mu(r)=\sqrt{\frac2{\pi}}\,\kappa^{-1} t^{-3/2}e^{-\kappa^2 t/2}\,
  e^{-r^2/2t}r\sinh\kappa r\,\d r$$
on $[0,\infty)$. By Lemma \ref{sect-4-lem-1}, $\mu$ is a
probability. Notice that the function $r\mapsto \log(1+2\kappa r)$
is concave, by Jensen's inequality,
  $$\eta_2(t)=\sqrt{\frac{\pi}2}\,\kappa t^{3/2}e^{\kappa^2 t/2}\int_0^\infty
  \log(1+2\kappa r)\,\d\mu(r)\leq \sqrt{\frac{\pi}2}\,\kappa t^{3/2}e^{\kappa^2 t/2}
  \log\bigg(1+2\kappa\int_0^\infty r\,\d\mu(r)\bigg).$$
Again by Lemma \ref{sect-4-lem-1},
  \begin{align*}
  \int_0^\infty r\,\d\mu(r)&=\sqrt{\frac2{\pi}}\,\kappa^{-1} t^{-3/2}e^{-\kappa^2 t/2}
  \int_0^\infty e^{-r^2/2t}r^2\sinh\kappa r\,\d r\cr
  &=\sqrt{\frac2{\pi}}\,t^{1/2}e^{-\kappa^2 t/2}+\sqrt{\frac2{\pi}}\,\kappa^{-1}(\kappa^2t+1)\alpha(t).
  \end{align*}
It is easy to show that $\sqrt{\frac2{\pi}}\,t^{1/2}
e^{-\kappa^2t/2}\leq \frac1\kappa \sqrt{\frac2{\pi e}}\leq
\frac1{2\kappa}$ and $\alpha(t)\leq \sqrt{\frac{\pi}2}$. Hence
  $$\int_0^\infty r\,\d\mu(r)\leq \frac1{2\kappa}+\kappa t+\frac1\kappa
  =\kappa t+\frac3{2\kappa}.$$
Consequently
  $$\eta_2(t)\leq \sqrt{\frac{\pi}2}\,\kappa t^{3/2}
  e^{\kappa^2 t/2}\log(2\kappa^2t+4).$$
Combining this with \eqref{sect-4-lem-3.1} and
\eqref{sect-4-lem-3.2}, we obtain the first inequality.

Now we prove the second inequality. By \eqref{sect-4.9}, we have
  \begin{equation}\label{sect-4-lem-3.3}
  \eta(t)<\eta_1(t)-\bar\eta_2(t),
  \end{equation}
where $\eta_1(t)$ is defined in \eqref{sect-4-lem-3.2} and
  $$\bar\eta_2(t)=\int_0^\infty e^{-r^2/2t}r\log(1+\kappa r)\sinh\kappa r\,\d r.$$
Define the measure
  $$\d\nu(r)=t^{-1/2}e^{-\kappa^2t/2}\alpha(t)^{-1}\, e^{-r^2/2t}\sinh\kappa r\,\d r$$
on $[0,\infty)$. Then by Lemma \ref{sect-4-lem-1}, $\nu$ is also a
probability. Notice that the function $r\mapsto r\log(1+\kappa r)$
is convex on $[0,\infty)$, again by Jensen's inequality,
  \begin{align*}
  \bar\eta_2(t)&=t^{1/2}e^{\kappa^2 t/2}\alpha(t)\int_0^\infty
  r\log(1+\kappa r)\,\d\nu(r)\cr
  &\geq t^{1/2}e^{\kappa^2 t/2}\alpha(t)\bigg(\int_0^\infty r\,\d\nu(r)\bigg)
  \log\bigg(1+\kappa \int_0^\infty r\,\d\nu(r)\bigg).
  \end{align*}
By Lemma \ref{sect-4-lem-1}, we have
  $$\int_0^\infty r\,\d\nu(r)=t^{-1/2}e^{-\kappa^2 t/2}\alpha(t)^{-1}
  \int_0^\infty e^{-r^2/2t}r\sinh\kappa r\,\d r
  =\sqrt{\frac{\pi}2}\,\kappa t\alpha(t)^{-1}.$$
Therefore
  $$\bar\eta_2(t)\geq \sqrt{\frac{\pi}2}\,\kappa t^{3/2}e^{\kappa^2t/2}
  \log\bigg(1+\sqrt{\frac{\pi}2}\,\kappa^2t\alpha(t)^{-1}\bigg).$$
Combining this with \eqref{sect-4-lem-3.3} and
\eqref{sect-4-lem-3.2} gives the second inequality.\fin

\medskip

Note that
  \begin{equation}\label{sect-4.12}
  \eta^\prime(t)=\frac1{2t^2}\int_0^\infty e^{-r^2/2t}r^3\sinh\kappa r
  \log\frac{\sinh\kappa r}{\kappa r}\,\d r.
  \end{equation}
In the same way we can prove the bounds on $\eta^\prime(t)$.

\begin{lemma}\label{sect-4-lem-4}
  \begin{align*}
  \eta^\prime(t)&> \frac12 \kappa^2t(\kappa^2t+5)+\frac12 \kappa t^{1/2}
  (\kappa^4t^2+6\kappa^2t+3)e^{\kappa^2 t/2}\alpha(t)\cr
  &\hskip12pt  - \frac12\sqrt{\frac{\pi}2}\,\kappa t^{1/2}(\kappa^2 t+3)e^{\kappa^2
  t/2}\cr
  &\hskip24pt\times \log\bigg(1+2\kappa\sqrt{\frac2{\pi}}\,\frac{t^{1/2}(\kappa^2t+5)}{\kappa^2t+3}e^{-\kappa^2t/2}
  +2\sqrt{\frac2{\pi}}\,\frac{\kappa^4t^2+6\kappa^2t+3}{\kappa^2t+3}\alpha(t)\bigg),\cr
  \eta^\prime(t)&< \frac12 \kappa^2t(\kappa^2t+5)+\frac12 \kappa t^{1/2}
  (\kappa^4t^2+6\kappa^2t+3)e^{\kappa^2 t/2}\alpha(t)\cr
  &\hskip12pt -\frac12\sqrt{\frac{\pi}2}\,\kappa t^{1/2}(\kappa^2t+3)e^{\kappa^2t/2}
  \log\bigg(1+\frac{\kappa^2t(\kappa^2t+3)e^{\kappa^2t/2}}
  {\kappa t^{1/2}+(\kappa^2t+1)e^{\kappa^2t/2}\alpha(t)}\sqrt{\frac{\pi}2}\,\bigg).
  \end{align*}
\end{lemma}

\noindent{\bf Proof.} The proofs are similar to Lemma
\ref{sect-4-lem-3}, hence we only give a sketch here. By
\eqref{sect-4.12} and \eqref{sect-4.9}, we have
  \begin{equation}\label{sect-4-lem-4.1}
  \eta^\prime(t)>\eta_3(t)-\eta_4(t),
  \end{equation}
where, by Lemma \ref{sect-4-lem-1},
  \begin{align}\label{sect-4-lem-4.2}
  \eta_3(t)&=\frac\kappa{2t^2}\int_0^\infty e^{-r^2/2t}r^4\sinh\kappa r\,\d
  r\cr
  &= \frac12 \kappa^2t(\kappa^2t+5)+\frac12 \kappa t^{1/2}
  (\kappa^4t^2+6\kappa^2t+3)e^{\kappa^2 t/2}\alpha(t)
  \end{align}
and
  $$\eta_4(t)=\frac1{2t^2}\int_0^\infty
  e^{-r^2/2t}r^3\log(1+2\kappa r)\sinh\kappa r\,\d r.$$
As in Lemma \ref{sect-4-lem-3}, using the concavity of $r\mapsto
\log(1+2\kappa r)$, we can get
  \begin{align*}
  \eta_4(t)&\leq \frac12\sqrt{\frac{\pi}2}\,\kappa t^{1/2}(\kappa^2 t+3)e^{\kappa^2
  t/2}\cr
  &\hskip12pt \times\log\bigg(1+2\kappa\sqrt{\frac2{\pi}}\,\frac{t^{1/2}(\kappa^2t+5)}{\kappa^2t+3}e^{-\kappa^2t/2}
  +2\sqrt{\frac2{\pi}}\,\frac{\kappa^4t^2+6\kappa^2t+3}{\kappa^2t+3}\alpha(t)\bigg).
  \end{align*}
Together with \eqref{sect-4-lem-4.1} and \eqref{sect-4-lem-4.2} gives us
the first inequality.

Next, by \eqref{sect-4.9},
  \begin{equation}\label{sect-4-lem-4.3}
  \eta^\prime(t)<\eta_3(t)-\bar\eta_4(t),
  \end{equation}
where $\eta_3(t)$ is defined in \eqref{sect-4-lem-4.2} and
  $$\bar\eta_4(t)=\frac1{2t^2}\int_0^\infty
  e^{-r^2/2t}r^3\log(1+\kappa r)\sinh\kappa r\,\d r.$$
Using the convexity of the function $r\mapsto r\log(1+\kappa r)$, we
can show that
  $$\bar\eta_4(t)\geq \frac12\sqrt{\frac{\pi}2}\,\kappa t^{1/2}(\kappa^2t+3)e^{\kappa^2t/2}
  \log\bigg(1+\frac{\kappa^2t(\kappa^2t+3)e^{\kappa^2t/2}}
  {\kappa t^{1/2}+(\kappa^2t+1)e^{\kappa^2t/2}\alpha(t)}\sqrt{\frac{\pi}2}\,\bigg).$$
Now the second inequality follows from the above estimate and
\eqref{sect-4-lem-4.3}, \eqref{sect-4-lem-4.2}. \fin

\medskip

Finally we are ready to prove Theorem \ref{sect-4-thm}.

\medskip

\noindent{\bf Proof of Theorem \ref{sect-4-thm}.} First we consider
the upper limit of the time derivative of the entropy. By
\eqref{sect-4.8} and Lemma \ref{sect-4-lem-3}, we have
  \begin{align*}
  \xi^\prime(t)\eta(t)< -\frac{\kappa}{\sqrt{2\pi}}\frac{\kappa^2t+3}{t^{1/2}e^{\kappa^2t/2}}
  -\frac{\alpha(t)}{\sqrt{2\pi}}\frac{\kappa^4t^2+4\kappa^2t+3}t
  +\frac{\kappa^2t+3}{2t}\log(2\kappa^2t+4).
  \end{align*}
Recall that $\alpha(t)$ is defined in Lemma \ref{sect-4-lem-1} and
$\lim_{t\ra\infty}\alpha(t)=\sqrt{\frac{\pi}2}$. And by Lemma
\ref{sect-4-lem-4}, we have
  \begin{align*}
  \xi(t)\eta^\prime(t)& < \frac{\kappa}{\sqrt{2\pi}}\frac{\kappa^2t+5}{t^{1/2}e^{\kappa^2t/2}}
  +\frac{\alpha(t)}{\sqrt{2\pi}}\frac{\kappa^4t^2+6\kappa^2t+3}t\cr
  &\hskip12pt -\frac{\kappa^2t+3}{2t}\log\bigg(1+\frac{\kappa^2t(\kappa^2t+3)e^{\kappa^2t/2}}{\kappa t^{1/2}
  +(\kappa^2t+1)e^{\kappa^2t/2}\alpha(t)}\sqrt{\frac{\pi}2}\,\bigg).
  \end{align*}
Summing up the above two estimates, we obtain
  \begin{align*}
  \frac{\d}{\d t}\big(\xi(t)\eta(t)\big)&<\sqrt{\frac2{\pi}}\,\frac{\kappa}{t^{1/2}e^{\kappa^2t/2}}
  +\sqrt{\frac2{\pi}}\,\kappa^2 \alpha(t)
  +\frac{\kappa^2t+3}{2t}
  \log\frac{2\kappa^2t+4}{1+\frac{\kappa^2t(\kappa^2t+3)e^{\kappa^2t/2}}{\kappa t^{1/2}
  +(\kappa^2t+1)e^{\kappa^2t/2}\alpha(t)}\sqrt{\frac{\pi}2}}.
  \end{align*}
From this it is easy to see that
  $$\varlimsup_{t\ra\infty}\frac{\d}{\d t}\big(\xi(t)\eta(t)\big)
  \leq \kappa^2+\frac{\kappa^2}2 \log 2=\kappa^2(1+\log \sqrt 2\,) .$$
Similarly we have
  $$\varliminf_{t\ra\infty}\frac{\d}{\d t}\big(\xi(t)\eta(t)\big)
  \geq \kappa^2-\frac{\kappa^2}2 \log 2=\kappa^2(1-\log \sqrt 2 \,) .$$
By \eqref{sect-4.2} and \eqref{sect-4.5}--\eqref{sect-4.7},
  $$\frac{\d}{\d t}\Ent(h(t,x,\cdot))=\frac 3{2t}+\kappa^2+\frac{\d}{\d t}\big(\xi(t)\eta(t)\big),$$
Noting that $\kappa^2=-k$, we finally complete the proof. \fin

\end{document}